\documentclass[11pt, a4paper]{article}
\usepackage{amsmath, amsthm, amssymb, eucal, setspace}

\makeatletter
\renewcommand\section{\@startsection{section}{1}{\z@}%
 						{-3.5ex \@plus -1ex \@minus -.2ex}
						{2ex \@plus.2ex}
						{\large\bfseries}}
\renewcommand\subsection{\@ifstar
						{\setcounter{subsection}{\value{equation}}
					\@startsection{subsection}{2}{\z@}
                          {1.75ex \@plus.5ex \@minus.2ex}%
                           {-.4em}		
					\textit*}
					{\setcounter{subsection}{\value{equation}}
						\stepcounter{equation}
					\@startsection{subsection}{2}{\z@}
                          {1.75ex \@plus.5ex \@minus.2ex}%
                           {-.4em}		
					\textit}}
\def\@seccntformat#1{\@ifundefined{#1@cntformat}%
	{\csname the#1\endcsname\quad} 
	{\csname #1@cntformat\endcsname}} 
\def\section@cntformat{\thesection.~} 
\def\subsection@cntformat{(\thesubsection)\ }
\renewcommand*\l@section{\mdseries\small\@dottedtocline{1}{1.5em}{2em}}
\makeatother

\numberwithin{equation}{section}
\theoremstyle{plain}
\newtheorem{maintheorem}{Theorem}

\swapnumbers

\theoremstyle{definition}

\theoremstyle{remark}
\newtheorem{remark}[equation]{Remark}

\newcommand{\cO}{\mathcal{O}}
\newcommand{\cE}{\mathcal{E}}
\newcommand{\cQ}{\mathcal{Q}}

\newcommand{\bP}{\mathbb{P}}

\newcommand{\bZ}{\mathbb{Z}}

\newcommand{\SL}{\mathrm{SL}}										
\newcommand{\SU}{\mathrm{SU}}
\newcommand{\GL}{\mathrm{GL}}
\newcommand{\Tr}{\mathrm{Tr}}
\newcommand{\Sym}{\mathrm{Sym}}
\newcommand{\ts}{\tilde\Sigma}										
										
\begin{document}
\title{\textbf{Branching of Hitchin's Prym cover for $\SL(2)$}}
\author{Constantin Teleman}
\date{\empty}
\maketitle

\begin{quote}
\abstract{\noindent
It is shown that the map from the Jacobian 
of the spectral curve to the moduli of stable bundles of rank $2$ 
is generically simply branched along an irreducible divisor. 
This observation falsifies the key step in the ``abelianization of 
the $\SU(2)$ WZW connection" presented in a recent paper \cite{yosh}. }
\end{quote}

\section{Statement}
Let $\Sigma$ be a smooth complex projective curve of genus $g\ge 3$ and 
$B$ a reduced divisor in $|K^2|$. The square root $r$ of a section of 
$K^2$ vanishing on $B$ defines a double cover $p:\ts\to\Sigma$ embedded 
in the total space of $K$, branched along $B$. It is a smooth curve of 
genus $\tilde{g} = 4g-3$, with Galois involution $\iota$, the sign 
change on $K$. $\ts$ is the simplest example of a \textit{spectral 
curve}  \cite{hit}, for rank $2$ bundles on $\Sigma$. More precisely, 
for a line bundle $L$ on $\ts$, the direct image $E=p_*L$ is a vector 
bundle on $\Sigma$, and multiplication by $r$ on sections of $L$ 
defines the \textit{Higgs field} $\phi: E\to E\otimes K$.

It is know that $E$  is stable, if $L$ avoids a sub-variety $V$ of 
co-dimension $\ge g-1$ in the Jacobian  of $\ts$ \cite{hit, bnr}. The 
construction works in families, so it defines a morphism $\pi$ from 
the Jacobian (minus $V$) to the moduli space of stable vector bundles 
on $\Sigma$.  Moreover, $\pi$ is generically finite, of degree $2^{3g-3}$. 
We chose $g\ge 3$ so that singularities of the moduli spaces, as well 
as the stable/semi-stable distinction can be ignored.

Let us concentrate on the critical Jacobian $\tilde J$ of degree 
$\tilde{g}-1$, which maps to the moduli space $M$ of semi-stable 
rank $2$ bundles of slope $g-1$; the story is similar for all even 
degrees. Call $K_M$ the canonical bundle of $M$. In this note, I 
verify the following (known) fact:

\begin{maintheorem} $\pi:\tilde{J}\setminus V\to M$ is \'etale  
away form an irreducible divisor $D$, and is generically simply 
branched along $D$. Moreover, $\cO(D) = \pi^*K_M^\vee$.
\end{maintheorem}
\noindent 
Up to isogeny, $\tilde{J}$ factors as $J\times P$ (see \ref{jacfactor}) 
and $D$ comes from an ample divisor on the Prym factor $P$. The 
important part is the simple branching; it implies the second statement, 
because the canonical bundle of $\tilde{J}$ is trivial and the Jacobian 
determinant of $\pi$ gives a section of  $\pi^*K_M^\vee$ with simple 
vanishing along $D$. 

In a recent paper \cite{yosh}, Yoshida proposed a solution of a 
long-standing problem, a reduction of the flat connection in the WZW 
model  for $\SU(2)$ to abelian Theta-functions. The key ingredient 
in the construction is a distinguished  Theta-fuction $\Pi$, living 
in a \textit{square root} of the anti-canonical pull-back $\pi^*K_M^\vee$ 
and vanishing along $D$. Both properties of $\Pi$ are essential for 
the constructions that follow. However, the theorem shows that 
such $\Pi$ does not exist.\footnote{Yoshida constructs $\Pi$ on an 
isogenous cover of $\tilde{J}$, but the distinction is unimportant. 
Page 2 of \textit{loc.~cit.} explicitly claims that $\Pi^2$ is the 
Jacobian determinant.} 
 
The interesting part of the story concerns $\SL(2)$ bundles and an 
associated Prym variety $P$; but their relation to $\GL(2)$ is 
straightforward, because $\pi$ is compatible with the tensor 
action (on $\tilde{J}$ and $M$) of the degree zero Jacobian $J$ of 
$\Sigma$. More precisely, let $\tilde{K}$ be the canonical bundle 
of $\ts$ and call $\tilde{B}$ the branch divisor; note the isomorphism 
$\tilde{K} \cong p^*K (\tilde{B})$. For a line bundle $L$ on 
$\ts$, the exact sequence 
\[
0\to L \to p^*p_*L \to \iota^*L (-\tilde{B}) \to 0
\]
shows the equivalence of the conditions 
\begin{equation}\label{prymcond}
L\otimes\iota^*L \cong \tilde{K} \quad\text{and }\quad 
	\det(p_*L) \cong K.
\end{equation}
They define the Prym variety $P\subset\tilde{J}$. Mind, however, that 
the first  isomorphism is always \textit{anti-invariant} for $\iota$, 
which changes the sign on the fibres of $\tilde{K}$ over $\tilde{B}$. 
With $M_K$ denoting the moduli space of semi-stable bundles on $\Sigma$ 
with determinant $K$ and $\Gamma\subset J$ its $2$-torsion subgroup, we 
have 
\begin{equation}\label{jacfactor}
\tilde{J} = J \times_\Gamma P \quad \text{and}  \quad 
	M = J\times_\Gamma M_K,
\end{equation}
compatibly with the map $\pi$. Up to translation, the restricted 
morphism $P\setminus V\to M_K$ is equivalent to the Prym covering 
of the moduli space of $\SL(2)$-bundles.  

\begin{remark}
$\tilde{K}\cong \cO(2\tilde{B})$, so one can use $L=\cO(\tilde{B})$ 
to identify $\tilde{J}$ with the degree zero Jacobian; $\iota^*$ 
becomes an automorphism.
\end{remark}

\section{Proof}
Let us abusively call the points in $\tilde{J}\setminus V$ where $\pi$ 
fails to be \'etale the `branch points', even though $\pi$ may not be 
everywhere finite; the contraction locus has co-dimension $\ge g-2$ 
(because the Theta-polarisations of the two spaces are compatible, 
Remark \ref{finalcomment}.i below). I describe the branching locus 
in terms of a ramified cover of a projective space and show its 
irreducibility. Finally, I show that the branching is simple by 
studying linearised deformations. 

\subsection{The branch locus.} 
Let us compare first-order deformations of $L$ and of $E=p_*L$. The 
tangent space to $P$ is the $(-1)$-eigenspace for $\iota$ on $H^1
(\ts;\cO)$, while the tangent space to $M_K$ at $E$ is $H^1(\Sigma; 
{\cE}nd^{\,0}(E))$, the traceless endomorphism bundle. 
Note that $p_*\cO$ splits into the $+/-$ eigenspaces of $\iota$ as 
$\cO\oplus K^\vee$, so that $TP$ is identified with $H^1(\Sigma; 
K^\vee)$. Unravelling the definition shows that the differential of 
$\pi$ at $L$ is the map induced by the Higgs field $\phi\in {\cE}nd^{\,0}
(E)\otimes K$:  
\[
\phi: H^1(\Sigma; K^\vee) \to H^1\left(\Sigma; {\cE}nd^{\,0}(E)\right).
\]
(For $\tilde{J}$ and $\GL(2)$, one adds the $H^1(\Sigma;\cO)$ summands 
to both sides.) When $E$ is stable, both spaces have the same dimension 
$3g-3$, and the short exact sequence on $\Sigma$,
\[
0\to K^\vee \xrightarrow{\phi} {\cE}nd^{\,0}(E) \to \cQ \to 0,
\]
shows that $\pi$ is not \'etale iff the quotient $\cQ$ has $h^1\neq 0$.  
In terms of $L$, $\cQ = p_*\left(\iota^*L^{-1}L(\tilde{B}) \right)$, and 
is a rank $2$ vector bundle with determinant $K$. It follows from Serre 
duality that $h^0(\cQ) = h^1(\cQ)$. Thus, $L$ is a branch point iff 
$\iota^*L^{-1}L(\tilde{B})$ has sections over $\ts$, in other words, 
the last line bundle lies in the Theta-divisor $\Theta$ of $\ts$.

\subsection{The Prym Theta-divisor.}
Consider the endomorphism $\sigma: L \mapsto \iota^*(L)^{-1}L(\tilde{B})$ 
of $\tilde{J}$. It factors via the projection to $\tilde{J}/J$ and 
lands in $P$. Restricted to $P$, $\sigma(L) = L^2(-\tilde{B})$ (or just 
the square, if we use $\cO(\tilde{B})$ as base-point). We now show 
that $\Theta$ meets $P$ transversely in an irreducible (and locally 
unibranch) divisor. Its pre-image $\sigma^*(\Theta\cap P)$ will be 
the branching divisor $D$ of $\pi$, and we will relate transversality 
to  simple branching.

Theta is the Abel-Jacobi image of $\Sym^{\tilde{g} -1}\ts$, and the 
condition $L\otimes\iota^*L \cong \tilde{K}$ defining $P$ says that 
each divisor $S\in |L|$ satisfies $S + \iota(S) \in |\tilde{K}|$: 
multiply the matching sections of $L$ and $\iota^*L$. The resulting 
section of $\tilde{K}$ is anti-invariant under $\iota$, as was the 
isomorphism in \eqref{prymcond}. The anti-invariant $p_*$-image 
of $\tilde{K}$ is $K^2$, and we obtain a bijection between divisors 
$S + \iota(S)\in |\tilde{K}|$ and points of $|K^2|$ (on $\Sigma$). 

Now, $S$ involves, in addition, a choice of point within each mirror 
pair in $S+\iota(S)$. The collection of choices defines a finite cover 
$\tilde{\bP}$ of $|K^2|$, simply branched over the hyperplanes of 
sections which vanish at some point of $B$. The monodromy around a 
hyperplane defined by $b\in B$ switches the point of $S$ which is 
near $b$ with its $\iota$-mirror. It follows that the monodromies 
act transitively on the fibres of $\tilde{\bP}\to |K^2|$, so that 
$\tilde{\bP}$ is irreducible. 
The same follows then for the intersection $\Theta\cap P$, which is 
set-theoretically the Abel-Jacobi image of $\tilde{\bP}$. Finally, 
the fibres of the Abel-Jacobi map are connected, so the image is 
locally unibranch. 

\subsection{Simple branching.} \label{simple}
First, observe that $P$ contains smooth points of $\Theta$. Indeed, 
over a singular point $L\in \Theta$, $\Sym^{\tilde{g} -1}\ts$ has 
positive-dimensional fibre; but this is also the fibre of the map 
$\tilde{\bP}\to\Theta\cap P$, which is generically finite for 
dimensional reasons. Next, at any smooth $L\in \Theta$ which lies 
in $P$, I claim that the normal to $\Theta$ is a $(-1)$-vector for 
$\iota$. For this, observe that the 
tangent space $T_L\Theta$ comprises the $\xi\in H^1(\ts;\cO)$ which 
induce the zero map $H^0(L) \to H^1(L)$, 
these $\xi$ being the first-order variations of $L$ which carry 
sections. Equivalently, the co-normal line to $\Theta$ is the image in 
$T^\vee\tilde{J} = H^0(\tilde{K})$ of the cup-product $H^0(L) \otimes 
H^0(\tilde{K}L^{-1})$. For $L\in P$, $\tilde{K}L^{-1} \cong \iota^* L$, 
so the image contains the product of a section with its $\iota$-transform; 
but we saw earlier that this is \textit {anti-invariant} under $\iota$. 
This proves transversality. 

In terms of $\pi$, this shows that $h^0(\Sigma; \cQ) = 1$ generically 
on $D$, and that the section fails to extend over the first-order 
neighbourhood of $D$ (which surjects to that of $\Theta\cap P$ in 
$P$). Since a first variation makes $\phi$ an isomorphism, the 
branching is simple. 

\subsection{Irreducibility.} Recall that in an Abelian variety of 
rank $2$ or more, any ample divisor is connected. As a connected 
\'etale cover of a locally unibranch divisor, $D$ is irreducible 
itself.  

\begin{remark}{\hspace{1cm}}\label{finalcomment}
\begin{enumerate}\itemsep0ex 
\item The moduli space $M$ is polarised by the inverse determinant of 
cohomology, which lifts to $\cO(\Theta)$ on $\tilde{J}$: this is because 
$H^*(\ts; L) = H^* (\Sigma; p_*L)$. However, $\cO(\Theta)$ is \textit
{not} principal on $P$. One way to normalise line bundles on $P$ is 
to relate them to $M_K$, whose Picard group is $\bZ$. The bundle 
$K_M^\vee$, which has Chern class $4$, lifts to $\sigma^*\cO(\Theta)$ 
over $P$. (This is the \textit{level\ $8$ line bundle} in 
\cite{yosh}.) 

\item
The sign in \S\ref{simple} is meaningful, as the opposite would make 
$\Theta$ tangent to $P$. Now, the Jacobian determinant of $\pi$ is the $\bar{\partial}$-determinant of $\cQ$. There is a perfect pairing 
$\cQ\otimes\cQ \to K$, the determinant; in terms of $\phi$,  
$q_1\wedge q_2 \mapsto \frac{1}{2}\Tr\left([\phi,q_1]\cdot 
q_2\right)$. The sign is in the skew-symmetry of the pairing; 
in the symmetric case, $\det\bar{\partial}$ would have a Pfaffian 
square root.
\end{enumerate}\end{remark}

\end{document}